\documentclass{ifacconf}
\pdfoutput=1

\counterwithin*{section}{part}

\usepackage[english]{babel}
\usepackage{amssymb}
\usepackage{tikz}
\usepackage{mathtools}
\usepackage{etoolbox}

\usepackage[ruled,vlined,algo2e,linesnumbered]{algorithm2e}
\makeatletter
\patchcmd{\algocf@Vline}{\vrule}{\vrule\vspace{-.32em}}{}{}
\makeatother
\SetStartEndCondition{ }{ }{}%
\SetKw{KwTo}{to}\SetKwFor{For}{for}{\string do}{}%
\SetKwIF{If}{ElseIf}{Else}{if}{then}{elif}{else}{}%
\SetKwFor{While}{while}{\string do}{}%

\DeclareMathAlphabet{\pazocal}{OMS}{zplm}{m}{n}
\usepackage [autostyle, english = american]{csquotes} 
\MakeOuterQuote{"}

\usepackage{xargs}
\usepackage[numbers]{natbib}
\usepackage{cleveref}
\usepackage{stmaryrd}

\usetikzlibrary{positioning,arrows,petri,calc,decorations.markings,arrows.meta}
\tikzset{
place/.style={circle,thick,minimum size=4mm,draw},
transitionV/.style={rectangle,thick,fill=black,minimum height=6mm,inner xsep=1pt}
}

\definecolor{myblue}{RGB}{0, 101, 202}
\definecolor{mygreen}{RGB}{130, 180, 0}
\definecolor{myred}{RGB}{197, 14, 31}
\definecolor{mypurple}{RGB}{128, 0, 128}
\definecolor{myyellow}{RGB}{204, 204, 0}
\definecolor{mygrey}{RGB}{105, 105, 105}

\crefname{thm}{theorem}{theorems}
\crefname{exmp}{example}{examples}

\newcommand{\dint}[1]{\left\llbracket#1\right\rrbracket} 

\newcommand{\K}{\mathcal{K}}
\renewcommand{\S}{\mathcal{S}}

\newcommand{\X}{\mathcal{X}}
\newcommand{\N}{\mathbb{N}}
\newcommand{\No}{\mathbb{N}_0}
\newcommand{\Z}{\mathbb{Z}}

\newcommand{\R}{\mathbb{R}}

\newcommand{\Rmax}{{\R}_{\normalfont\fontsize{7pt}{11pt}\selectfont\mbox{max}}}

\newcommand{\Rbar}{\overline{\R}}

\newcommand{\rewrite}[2]{} 

\newcommand{\graph}{\mathcal{G}}

\DeclareMathOperator{\Ima}{Im}

\makeatletter
\newcommand{\splus}{%
  \DOTSB\mathop{\mathpalette\mattos@splus\relax}\slimits@
}
\newcommand\mattos@splus[2]{%
  \vcenter{\hbox{%
    \sbox\z@{$#1\oplus$}%
    \resizebox{!}{0.9\dimexpr\ht\z@+\dp\z@}{\raisebox{\depth}{$\m@th#1\boxplus$}}%
  }}%
  \vphantom{\oplus}%
}
\makeatother

\makeatletter
\newcommand{\stimes}{%
  \DOTSB\mathop{\mathpalette\mattos@stimes\relax}\slimits@
}
\newcommand\mattos@stimes[2]{%
  \vcenter{\hbox{%
    \sbox\z@{$#1\otimes$}%
    \resizebox{!}{0.9\dimexpr\ht\z@+\dp\z@}{\raisebox{\depth}{$\m@th#1\boxtimes$}}%
  }}%
  \vphantom{\otimes}%
}
\makeatother

\usepackage{pict2e}

\makeatletter
\newcommand*{\bigsplus}{\DOTSB\mathop{\mathpalette\big@boxplus\relax}\slimits@}

\newcommand{\big@boxplus}[2]{%
  \vcenter{%
    \m@th\bigbox@thickness{#1}%
    \sbox\z@{$#1\bigoplus$}%
    \dimen@=\ht\z@ \advance\dimen@\dp\z@
    \hbox{%
      \setlength{\unitlength}{\dimen@}%
      \begin{picture}(1,1)
      \polyline(0.1,0.1)(0.9,0.1)(0.9,0.9)(0.1,0.9)(0.1,0.1)(0.5,0.1)
      \polyline(0.5,0.1)(0.5,0.9)
      \polyline(0.1,0.5)(0.9,0.5)
      \end{picture}%
    }%
  }%
}

\newcommand{\bigbox@thickness}[1]{%
  \ifx#1\displaystyle
    \linethickness{0.2ex}%
  \else
    \ifx#1\textstyle
      \linethickness{0.16ex}%
    \else
      \ifx#1\scriptstyle
        \linethickness{0.12ex}%
      \else
        \linethickness{0.1ex}%
      \fi
    \fi
  \fi
}
\makeatother

\newcommand{\vsdots}{\raisebox{3pt}{$\scalebox{.75}{\vdots}$}}

\renewcommand{\epsilon}{\varepsilon}

\newsavebox{\rightdiv}
\sbox\rightdiv{\tikz[anchor=south,baseline]{\footnotesize\node[inner sep=0pt,minimum height=.75em] at (0,0) (A){$\circ$};\node[inner sep=0pt,minimum height=.75em] at (0,0) {$/$};}}
\newsavebox{\rightmindiv}
\sbox\rightmindiv{\tikz[anchor=south,baseline]{\footnotesize\node[inner sep=0pt,minimum height=.75em] at (0,0) (A){$\bullet$};\node[inner sep=0pt,minimum height=.75em] at (0,0) {$/$};}}

\newsavebox{\leftdiv}
\sbox\leftdiv{\tikz[anchor=south,baseline]{\footnotesize\node[inner sep=0pt,minimum height=.75em] at (0,0) (A){$\circ$};\node[inner sep=0pt,minimum height=.75em] at (0,0) {$\setminus$};}}
\newsavebox{\leftmindiv}
\sbox\leftmindiv{\tikz[anchor=south,baseline]{\footnotesize\node[inner sep=0pt,minimum height=.75em] at (0,0) (A){$\bullet$};\node[inner sep=0pt,minimum height=.75em] at (0,0) {$\setminus$};}}

\newcommand{\transposed}{^\intercal}

\newcommand{\ie}{i.e.,\ }

\newcommand{\qedhere}{\tag*{$\blacksquare$}} 

\pgfdeclarelayer{back}
\pgfsetlayers{back,main}
\makeatletter
\pgfkeys{%
  /tikz/on layer/.code={
    \pgfonlayer{#1}\begingroup
    \aftergroup\endpgfonlayer
    \aftergroup\endgroup
  },
  /tikz/node on layer/.code={
    \pgfonlayer{#1}\begingroup
    \expandafter\def\expandafter\tikz@node@finish\expandafter{\expandafter\endgroup\expandafter\endpgfonlayer\tikz@node@finish}%
  },
}
\tikzset{%
glow/.style={%
preaction={#1, draw, line cap=round, line join=round, line width=0.5pt, opacity=0.04, on layer=back,
preaction={#1, draw, line cap=round, line join=round, line width=1.0pt, opacity=0.04, on layer=back,
preaction={#1, draw, line cap=round, line join=round, line width=1.5pt, opacity=0.04, on layer=back,
preaction={#1, draw, line cap=round, line join=round, line width=2.0pt, opacity=0.04, on layer=back,
preaction={#1, draw, line cap=round, line join=round, line width=2.5pt, opacity=0.04, on layer=back,
preaction={#1, draw, line cap=round, line join=round, line width=3.0pt, opacity=0.04, on layer=back,
preaction={#1, draw, line cap=round, line join=round, line width=3.5pt, opacity=0.04, on layer=back,
preaction={#1, draw, line cap=round, line join=round, line width=4.0pt, opacity=0.04, on layer=back,
preaction={#1, draw, line cap=round, line join=round, line width=4.5pt, opacity=0.04, on layer=back,
preaction={#1, draw, line cap=round, line join=round, line width=5.0pt, opacity=0.04, on layer=back,
preaction={#1, draw, line cap=round, line join=round, line width=5.5pt, opacity=0.04, on layer=back,
preaction={#1, draw, line cap=round, line join=round, line width=6.0pt, opacity=0.04, on layer=back,
}}}}}}}}}}}}}}

\begin{document}
 

\begin{frontmatter}

    \title{Controlled Invariance\\in Fully Actuated Max-plus Linear Systems\\with Precedence Semimodules\thanksref{footnoteinfo}}

\thanks[footnoteinfo]{Support from Deutsche Forschungsgemeinschaft (DFG) via grant RA 516/14-1 and under Germany’s Excellence Strategy -- EXC 2002/1 ``Science of Intelligence'' -- project number 390523135 is gratefully acknowledged.}

\author[1]{Davide Zorzenon}
\author[1,2]{J\"{o}rg Raisch}

\address[1]{Control Systems Group, Technische Universit\"at Berlin, Germany (e-mail: [zorzenon,raisch]@control.tu-berlin.de)}
\address[2]{Science of Intelligence, Research Cluster of Excellence, Berlin, Germany}

\begin{abstract}
Given a max-plus linear system and a semimodule, the problem of computing the maximal controlled invariant subsemimodule is still open to this day.
In this paper, we consider this problem for the specific class of fully actuated systems and constraints in the form of precedence semimodules.
The assumption of full actuation corresponds to the existence of an input for each component of the system state.
A precedence semimodule is the set of solutions of inequalities typically used to represent time-window constraints.
We prove that, in this setting, it is possible to (i) compute the maximal controlled invariant subsemimodule and (ii) decide the convergence of a fixed-point algorithm introduced by R.D. Katz in strongly polynomial time.
\end{abstract}

\begin{keyword}
Max-plus algebra, controlled invariance, graph theory, P-time event graphs
\end{keyword}

\end{frontmatter}

\section{Introduction}

The concept of controlled invariance constitutes the cornerstone of the so-called geometric approach in control theory.
Given a dynamical system, we say that a set -- representing specifications for the system state -- is controlled invariant if there exists a control action that keeps the state inside this set.
When the considered dynamical system is linear over a field and the set forms a vector space, there exists an efficient fixed-point algorithm that computes the maximal controlled invariant subspace.
This result is the basic ingredient for the solution to numerous control design problems, such as the disturbance decoupling problem and the model matching problem (see \cite{basile2002controlled,wonham1974linear}).

In linear systems over the max-plus semiring (or max-plus linear systems), the situation is more involved because finite-time convergence of the same fixed-point algorithm mentioned above, adjusted to the different algebraic setting, is not guaranteed (see \cite{katz2007maxplus}).
Due to this difficulty, reminiscent of the case over rings (\cite{conte1995disturbance}), the complexity of solving several control problems in manufacturing and transportation networks modeled as max-plus linear systems is, as of today, unknown.
Among the applications of geometric methods in this framework we mention the design of dynamic observers (\cite{di2010duality}), the enforcement of time-window constraints (\cite{maia2011control}) and more generic constraints in the steady state (\cite{goncalves2016onthe}), and the solution of the model matching problem (\cite{martinez2022systems,9832498}).

This paper presents an interesting class of max-plus linear systems and semimodules\footnote{Semimodules over semirings are defined analogously to vector spaces over fields.} for which the convergence of the fixed-point algorithm can be verified in strongly polynomial time complexity.
Specifically, we consider \emph{fully actuated} max-plus linear systems (\ie in which the input matrix of the state-space model is the identity, see \eqref{eq:fully_actuated}) where the trajectory of the state needs to satisfy \emph{precedence constraints} (\ie inequalities of the form $x_i \geq A_{ij} + x_j$, see \eqref{eq:linear_dual_inequalities}).
In typical applications of such systems, precedence constraints represent time-window constraints (as those considered in the article \cite{maia2011control}) and the assumption of full actuation translates into the ability to delay at will the occurrence of every event in the system.
Although the latter assumption is clearly restrictive, we observe that this condition is satisfied in applications such as robotic job-shops and transportation systems (an example is given in \Cref{se:5}). 

After recalling in \Cref{se:2} the basics of max-plus algebra and precedence constraints, \Cref{se:3} shows how to decide, for the considered class of systems, whether a control action exists under which the state satisfies all constraints.
The strategies presented in \Cref{se:3} are based on recent results on the analysis of consistency in P-time event graphs (see \cite{ZORZENON202219,zorzenon2024infinite}).
In \Cref{se:4}, we interpret these discoveries in the framework of the geometric approach.
The results of \Cref{se:4} show a deep connection between the concept of controlled invariance and the longest path problem in infinite graphs, which hopefully will bring useful insights for the solution of similar problems in more complex settings.
\Cref{se:5} presents simple applications examples (including a railway network adjusted from \cite{katz2007maxplus}) that are not solvable using previous techniques, and \Cref{se:conclusions} gives suggestions for future work.

\subsection*{Notation}

$\Z$ and $\R$ are the sets of integers and reals, respectively.
The sets of positive and non-negative integers are denoted respectively by $\N$ and $\No$.
Given two integers $a,b\in\Z$ such that $a\leq b$, $\dint{a,b}$ indicates the set $\{a,a+1,a+2,\dots,b\}$.
Moreover, $\Rmax \coloneqq \R \cup \{-\infty\}$, and $\Rbar \coloneqq \R\cup \{-\infty,+\infty\}$.

\section{Preliminaries}\label{se:2}

\subsection{Max-plus algebra}

The max-plus algebra is the set $\Rbar$ endowed with operations $\oplus$ and $\otimes$, defined by: for all $a,b\in\Rbar$,
\[
    a\oplus b = \max\{a,b\},\quad a\otimes b = 
    \begin{dcases}
        a+b & \mbox{if $a,b\in\R\cup\{+\infty\}$},\\
        -\infty & \mbox{if $a$ or $b$ is $-\infty$.}
    \end{dcases}
\]
Given two (possibly infinite) sets $M,N\subseteq\N$, a max-plus matrix with index sets $M$ and $N$ is a function $A:M\times N\rightarrow \Rbar$.
We denote by $\Rbar^{M\times N}$ the set of all max-plus matrices with index sets $M$ and $N$, and by $A_{ij}$ the value $A(i,j)$.
Similarly, $\Rmax^{M\times N}$, respectively $\R^{M\times N}$, are the sets of all matrices with index sets $M$ and $N$ and entries in $\Rmax$, respectively $\R$.
If $N$ is a singleton, we write $\Rbar^M$, respectively $\Rmax^M$ or $\R^M$.
If $M=\dint{1,m}$ and $N=\dint{1,n}$, we will write $\Rbar^{m\times n}$ instead of $\Rbar^{M\times N}$.
Operations $\oplus$ and $\otimes$ are extended to matrices in the usual way: given $A,B\in\Rbar^{M\times N}$ and $C\in\Rbar^{N\times P}$,
\[
    (A\oplus B)_{ij} = A_{ij}\oplus B_{ij},\quad (A\otimes C)_{ij} = \bigoplus_{k\in N} A_{ik} \otimes C_{kj}.
\]
In $\Rbar^{N\times N}$, the neutral elements for operations $\oplus$ and $\otimes$ are, respectively, the zero matrix $\mathcal{E}$ and the identity matrix $E$, defined such that $\mathcal{E}_{ij} = -\infty$ for all $i,j$, $E_{ij} = 0$ if $i=j$ and $E_{ij} = -\infty$ if $i\neq j$.
Sometimes, we will use the notation $\mathcal{E}_{m\times n}$ to indicate a matrix of size $m\times n$ containing only $-\infty$'s.
Given a scalar $\lambda\in\Rbar$ and a matrix $A\in\Rbar^{M\times N}$, the scalar-matrix multiplication $\lambda\otimes A$ results in a matrix with index sets $M$ and $N$ and coefficients $(\lambda\otimes A)_{ij} = \lambda\otimes A_{ij}$.
To simplify notation, we will often omit the symbol $\otimes$ and write $a b$ in place of $a\otimes b$.
Given two matrices $A,B\in\Rbar^{M\times N}$, we write $A\leq B$ to indicate that, for all $i,j$, $A_{ij} \leq B_{ij}$.
The Kleene star $A^*$ of matrix $A\in\Rbar^{N\times N}$ is defined by
\[
    A^* = \bigoplus_{k=0}^{+\infty} A^k,
\]
where $A^0 = E$ and $A^{k} = A^{k-1}\otimes A$ for all $k\in\N$.

A subset $\X$ of $\Rmax^n$ that is closed under finitely many additions and multiplications by scalars from $\Rmax$ is called semimodule.
A semimodule $\X\subseteq \Rmax^n$ is finitely generated if there exists a matrix $U\in\Rbar^{n\times p}$ such that $\X$ is the image of $U$, \ie $\X = \Ima U \coloneqq \{ U\otimes u \mid u\in\Rmax^p \}$.
For all matrices $A,B\in\Rmax^{N\times N}$, one has $(A^*)^* = A^*A^* = A^*$ and, if $A^* \in \Rmax^{N\times N}$, then $\Ima A^* = \Ima B^* \ \Leftrightarrow\ A^* = B^*$.
A square matrix $A$ such that $A^* = A$ is called a \emph{star matrix}.

\subsection{Precedence constraints and graphs}\label{su:precedence_constraints}

Precedence constraints are systems of (finitely or infinitely many) inequalities of the form, for all $i,j\in N$, $x_i \geq A_{ij} + x_j$, where $A\in\Rmax^{N\times N}$ and $x\in\R^{N}$.
In the max-plus algebra, precedence constraints can be written as
\begin{equation}\label{eq:precedence_constraint_max_plus}
    x \geq A\otimes x.
\end{equation}
From \cite[Equation 6.11]{hardouin2018control}, the latter inequality has two equivalent expressions: $x\geq A^* \otimes x$ and $x = A^* \otimes x$.
Moreover, $x$ is a solution if and only if $x$ belongs to the \emph{precedence semimodule} $\Ima A^*$.\footnote{Indeed, $x\in\Ima A^*$ implies $x = A^* u$ for some $u$ and, because $A^* = A^* A^*$, $x = A^* A^* u = A^* x$. On the other hand, if $x = A^* x$, then clearly $x\in\Ima A^*$.}

It is convenient to represent precedence constraints by means of precedence graphs.
The precedence graph $\graph(A)$ associated with \eqref{eq:precedence_constraint_max_plus} is a weighted directed graph with a node for each element of $N$ and an arc from node $j$ to node $i$ of weight $A_{ij}$ if and only if $A_{ij}\neq -\infty$.
Recall that $(A^\ell)_{ij}$ is equal to the supremal weight of all paths $\rho$ in $\graph(A)$ from node $j$ to node $i$ of length $\ell$.
Therefore, if $\rho$ is the path with maximum weight $|\rho|_W$ among all the paths from node $j$ to node $i$, then $|\rho|_W = (A^{*})_{ij}$.

The following theorem contains a fundamental observation by \cite{gallai1958maximum}, generalized to the case of infinitely many precedence constraints in \cite{zorzenon2024infinite}.
\begin{thm}\label{th:gallai}
The precedence constraints \eqref{eq:precedence_constraint_max_plus} admit a real solution $x\in\R^N$ if and only if the Kleene star $A^*$ of matrix $A$ converges in $\Rmax^{N\times N}$, \ie $(A^*)_{ij}\in\Rmax$ for all $i,j\in N$.
\end{thm}
In terms of the precedence graph $\graph(A)$, this condition is equivalent to the absence of $\infty$-weight paths, \ie sequences $\rho_1,\rho_2,\ldots$ of paths connecting two nodes of $\graph(A)$ with infinite limit weight, $\lim_{k\rightarrow +\infty}|\rho_k|_W = +\infty$. 

In finite precedence graphs, the presence of an $\infty$-weight path is a necessary and sufficient condition for the existence of a positive-weight circuit, \ie of a path $\rho$ with weight $|\rho|_W>0$ starting and ending at the same node $i$. 
For instance, the circuit $\rho_1 = 1\rightarrow 2 \rightarrow 1$ in \Cref{fi:finite_precedence_graph} has positive weight $|\rho_1|_W = 1$, which implies that the sequence of paths $\rho_1,\rho_2,\ldots$, where $\rho_k$ is obtained by repeating $k$ times the circuit $\rho_1$, has infinite limit weight.
Recall that the existence of positive-weight circuits in precedence graphs with $n$ nodes can be verified in strongly polynomial time $O(n^3)$ using, for instance, the Floyd-Warshall algorithm (\cite{cormen2022introduction}).

\begin{figure}[h]
    \centering
    \begin{tikzpicture}[node distance=1cm and 1cm,>=stealth',bend angle=45,double distance=.5mm,arc/.style={->,>=stealth'},place/.style={circle,thick,fill=gray!15,minimum size=4mm,draw}]

\tiny

\node [place] (a) {$1$};
\node [place,right=of a] (b) {$2$};

\draw [arc] (a) to [bend left=30] node [auto] {$2$} (b);
\draw [arc] (a) to [loop,out=180-22.5,in=180+22.5,looseness=6] node [left] {$-3$} (a);
\draw [arc] (b) to [bend left=30] node [above] {$-1$} (a);

\end{tikzpicture}
    \caption{Finite precedence graph with an $\infty$-weight path.}\label{fi:finite_precedence_graph}
\end{figure}

In infinite precedence graphs, there are other ways to generate an $\infty$-weight path.
For example, the precedence graph in \Cref{fi:infinite_precedence_graph} contains an $\infty$-weight path from node $1$ to node $2$, but no positive-weight circuit.
A path $\rho_k$ from the sequence $\rho_1,\rho_2,\ldots$ with infinite limit weight can be defined by $\rho_k = 1\rightarrow 3 \rightarrow 5 \rightarrow \ldots \rightarrow 2k-1 \rightarrow 2k \rightarrow 2k-2 \rightarrow \ldots \rightarrow 4 \rightarrow 2$; since $|\rho_k|_W = k-1$, we have $\lim_{k\rightarrow+\infty}|\rho_k|_W = +\infty$.

\begin{figure}[ht]
    \centering
    \begin{tikzpicture}[node distance=2cm and 2cm,>=stealth',bend angle=45,double distance=.5mm,arc/.style={->,>=stealth'},place/.style={circle,thick,minimum size=4mm,draw,fill=gray!15}]

\tiny

\foreach \z in {1,2,3,4,5}
{
\pgfmathtruncatemacro{\za}{\z*2-1}
\pgfmathtruncatemacro{\zb}{\z*2}
\node [place] (ntop\z) at (\z,0) {$\za$};
\node [place] (nbot\z) at (\z,-1) {$\zb$};
\draw [arc] (ntop\z) to node[auto] {$0$} (nbot\z);
}
\draw [arc] (ntop1) to node[auto] {$2$} (ntop2);
\draw [arc] (ntop2) to node[auto] {$2$} (ntop3);
\draw [arc] (ntop3) to node[auto] {$2$} (ntop4);
\draw [arc] (ntop4) to node[auto] {} (ntop5);

\draw [arc] (nbot2) to node[auto] {$-1$} (nbot1);
\draw [arc] (nbot3) to node[auto] {$-1$} (nbot2);
\draw [arc] (nbot4) to node[auto] {$-1$} (nbot3);
\draw [arc] (nbot5) to node[auto] {} (nbot4);

\fill [white] (4.5,.5) rectangle (5.5,-1.5);
\node (dots2) at (5,-.5) {\textbf{\dots}};

\draw [arc] (1,-1.5) to node[pos=.95,above] {$\N$} (5.5,-1.5);
\end{tikzpicture}
    \caption{Infinite precedence graph with an $\infty$-weight path but no positive-weight circuit.}\label{fi:infinite_precedence_graph}
\end{figure}

\section{Problem statement and solution}\label{se:3}

\subsection{The problem}

A max-plus linear system is a dynamical system evolving according to
\begin{equation}\label{eq:max_plus_linear_system}
    (\forall k\in\N)\quad x(k+1) = A x(k) \oplus B u(k),
\end{equation}
where $x:\N\rightarrow \Rmax^n$, $A\in\Rmax^{n\times n}$, $B\in\Rmax^{n\times m}$, and $u:\N\rightarrow \Rmax^m$.
In the context of discrete-event systems, typically $x_i(k)$ has the meaning of "time instant of the $k$-th occurrence of event $i$".

Let $A,L,C,\tilde{R}$ be four $n\times n$ matrices with elements from $\Rmax$.
Consider a fully actuated max-plus linear system, \ie a system \eqref{eq:max_plus_linear_system} in which matrix $B$ is the max-plus identity matrix of dimension $n$, \ie $B = E\in\Rmax^{n\times n}$:
\begin{equation}\label{eq:fully_actuated}
    (\forall k\in\N)\quad x(k+1) = Ax(k) \oplus u(k).
\end{equation}
We want to impose that the trajectory $\{x(k)\}_{k\in\N}$ satisfies the following inequalities:
\begin{equation}\label{eq:linear_dual_inequalities}
    (\forall k\in\N)\quad
    \left\{
\begin{array}{rcl}
      x(k) &\geq& L\otimes x(k+1),\\
      x(k) &\geq& C\otimes   x(k),\\
    x(k+1) &\geq& \tilde{R}\otimes   x(k).
\end{array}
    \right.
\end{equation}
With simple manipulations, the above inequalities can be written in the standard algebra as
\begin{equation*}
    \left(
\begin{array}{c}
    \forall k\in\N,\\
\forall i,j\in\dint{1,n}
\end{array} 
    \right)\
    \left\{
\begin{array}{rcl}
    C_{ij}  \leq & x_i(k) - x_j(k)& \leq -C_{ji},\\
    \tilde{R}_{ij}  \leq & x_i(k+1) - x_j(k) &\leq -L_{ji}.
\end{array}
    \right.
\end{equation*}
This formulation shows more clearly that the inequalities in \eqref{eq:linear_dual_inequalities} can be interpreted as time-window constraints on the occurrence of events in a discrete-event system.

The main decision problem that we consider in this paper is to determine, for given matrices $A,L,C,\tilde{R}$, if there exists a sequence $\{u(k)\}_{k\in\N}$ such that $\{x(k)\}_{k\in\N}$ satisfies all inequalities \eqref{eq:linear_dual_inequalities}.
Based on recent results on the analysis of P-time event graphs, we will show that this problem can be solved in strongly polynomial time.

\subsection{Equivalent formulation}\label{su:equivalent_formulation}

Note that, given a vector $x(k)$, the set of possible $x(k+1)$ according to the expression \eqref{eq:fully_actuated} can be re-stated as
\begin{multline*}
             \{x(k+1) \in\Rmax^n\mid (\exists u(k)\geq Ax(k))\ x(k+1) = u(k)\}.
\end{multline*}
This is because $x(k+1)\geq Ax(k)$ and we can assume without loss of generality, since $u(k)$ is free, that $u(k)\geq Ax(k)$, from which $x(k+1) = u(k)$ follows immediately.
Summarizing, for $u(k)\geq Ax(k)$, expressions \eqref{eq:fully_actuated} and \eqref{eq:linear_dual_inequalities} are equivalent to
\begin{equation}\label{eq:scenario_geometric_approach}
    (\forall k\in\N)\quad
    \left\{
    \begin{array}{rcl}
        x(k+1) &=& u(k),\\
      x(k+1) &\geq& A\otimes x(k),\\
      x(k) &\geq& L\otimes x(k+1),\\
      x(k) &\geq& C\otimes   x(k),\\
    x(k+1) &\geq& \tilde{R}\otimes   x(k).
    \end{array} 
    \right.
\end{equation}
Note that the conjunction of $x(k+1)\geq A x(k)$ and $x(k+1)\geq \tilde{R} x(k)$ is equivalent to the inequality $x(k+1) \geq A x(k) \oplus \tilde{R} x(k) = R x(k)$, where $R\coloneqq A\oplus \tilde{R}$.
Moreover, observe that $u(k)$ plays no role in \eqref{eq:scenario_geometric_approach}, and thus the equation $x(k+1)=u(k)$ can be eliminated. 
In conclusion, \eqref{eq:fully_actuated} and \eqref{eq:linear_dual_inequalities} can be restated as
\begin{equation}\label{eq:P-TEG}
    (\forall k\in\N)\quad
    \left\{
    \begin{array}{rcl}
      x(k) &\geq& L\otimes x(k+1),\\
      x(k) &\geq& C\otimes   x(k),\\
    x(k+1) &\geq& R\otimes   x(k).
    \end{array} 
    \right.
\end{equation}

\subsection{The solution}\label{su:solution}

The system of inequalities \eqref{eq:P-TEG} represents the dynamics of P-time event graphs, a class of time discrete-event systems introduced in \cite{khansa1996a}.
Here we collect the main results related to the existence of trajectories $\{x(k)\}_{k\in\N}$ satisfying \eqref{eq:P-TEG}, from which the solution of our problem follows immediately.

We say that system \eqref{eq:P-TEG} is \emph{consistent} if it admits an infinite trajectory $\{x(k)\}_{k\in\N}$ that satisfies all inequalities.
We also define a weaker property, called \emph{weak consistency}, which will be useful in \Cref{se:4}.
System \eqref{eq:P-TEG} is weakly consistent if, for all $K\in\N$, there exists a finite trajectory $\{x(k)\}_{k\in\dint{1,K}}$ satisfying
\begin{equation}\label{eq:P-TEG_finite}
    \begin{array}{rrcl}
        (\forall k\in\dint{1,K-1}) \quad&
      x(k) &\geq& L\otimes x(k+1),\\
        (\forall k\in\dint{1,K})  \quad&
      x(k) &\geq& C\otimes   x(k),\\
        (\forall k\in\dint{1,K-1})  \quad&
    x(k+1) &\geq& R\otimes   x(k).
    \end{array} 
\end{equation}

System \eqref{eq:P-TEG_finite} can be rewritten as the precedence inequality $x_{\dint{K}} \geq M_{\dint{K}} \otimes x_{\dint{K}}$, where $M_{\dint{K}}\in\Rmax^{Kn\times Kn}$ and $x_{\dint{K}}\in\R^{Kn}$ are defined by
\[
    M_{\dint{K}} = 
    \begin{bmatrix}
        C & L & \mathcal{E} & \mathcal{E} & \cdots & \mathcal{E} &  \mathcal{E}\\
        R & C & L & \mathcal{E} & \cdots & \mathcal{E} &  \mathcal{E}\\
        \mathcal{E} & R & C & L & \cdots & \mathcal{E} &  \mathcal{E}\\
        \mathcal{E} & \mathcal{E} & R & C & \cdots & \mathcal{E} &  \mathcal{E}\\
        \vdots & \vdots & \vdots & \vdots & \ddots & \vdots & \vdots\\
        \mathcal{E} & \mathcal{E} & \mathcal{E} & \mathcal{E} & \cdots & C &  L\\
        \mathcal{E} & \mathcal{E} & \mathcal{E} & \mathcal{E} & \cdots & R &  C
    \end{bmatrix},\
    x_{\dint{K}} = 
    \begin{bmatrix}
    x(1)\\x(2)\\x(3)\\x(4)\\\vdots\\x(K-1) \\ x(K)
    \end{bmatrix}.
\]
Therefore, weak consistency is equivalent to the absence of positive-weight circuits in $\graph(M_{{\dint{K}}})$, for all $K\in\N$.
Similarly, \eqref{eq:P-TEG} is equivalent to the precedence inequality $x_{\dint{\infty}}\geq M_{\dint{\infty}} \otimes x_{\dint{\infty}}$, where $M_{\dint{\infty}}\in\Rmax^{\N\times \N}$ and $x_{\dint{\infty}}\in\R^{\N}$ are defined as the limit, for $K\rightarrow +\infty$, of $M_{\dint{K}}$ and $x_{\dint{K}}$, respectively, \ie
\[
    M_{\dint{\infty}} = 
    \begin{bmatrix}
        C & L & \mathcal{E} & \mathcal{E} & \cdots \\ 
        R & C & L & \mathcal{E} & \cdots \\
        \mathcal{E} & R & C & L & \cdots \\
        \mathcal{E} & \mathcal{E} & R & C & \cdots \\
        \vdots & \vdots & \vdots & \vdots & \ddots 
    \end{bmatrix},\ 
    x_{\dint{\infty}}= 
   \begin{bmatrix}
       x(1)\\x(2)\\x(3)\\x(4)\\\vdots
   \end{bmatrix}.
\]

Let $\Pi_0,\Pi_1,\ldots$ 
be the sequence of matrices in $\Rbar^{n\times n}$ defined recursively by
\[
    \Pi_0 = C^*,\quad \Pi_{k+1} = (L \Pi_k R\oplus C)^*.
\]
Moreover, let $\Pi_\infty=\lim_{k\rightarrow\infty} \Pi_k$.
We recall the following result, obtained in \cite{ZORZENON202219,zorzenon2024infinite}.

\begin{thm}\label{th:consistency}
    System \eqref{eq:P-TEG} is consistent if and only if $\graph(M_{\dint{\infty}})$ does not contain $\infty$-weight paths, and is weakly consistent if and only if $\graph(M_{\dint{\infty}})$ does not contain positive-weight circuits.
Consistency can be checked in time $O(n^5)$ and weak consistency in $O(n^9)$.
In particular,
\begin{enumerate}
    \item consistency is equivalent to the conditions: $\Pi_{n^2+1} = \Pi_{n^2} \mbox{ and } \Pi_{n^2}\in\Rmax^{n\times n}$,
    \item \eqref{eq:P-TEG} is weakly consistent but not consistent if and only if, for all $k\in\No$, $\Pi_k\in\Rmax^{n\times n}$ and $\Pi_{k+1}\neq \Pi_{k}$,
    \item weak consistency is equivalent to the condition: for all $k\in\N$, $\Pi_{k}\in\Rmax^{n\times n}$.
\end{enumerate}
\end{thm}

\begin{rem}
    Note that if $\Pi_{n^2+1} = \Pi_{n^2}$, then $\Pi_{k} = \Pi_{n^2}$ for all $k\geq n^2$, including $k=\infty$.
\end{rem}

\begin{exmp}\label{ex:cont}
    Consider the fully actuated max-plus linear system \eqref{eq:fully_actuated} subject to constraints \eqref{eq:linear_dual_inequalities} with matrices 
    \[
        A = \begin{bmatrix}
            2 & -\infty\\
            -\infty & -\infty
        \end{bmatrix},\ 
        L = \begin{bmatrix}
            -\infty & -\infty\\
            -\infty & -1
        \end{bmatrix},\ 
        C = \begin{bmatrix}
            -\infty & -\infty\\
            0 & -\infty\\
        \end{bmatrix},
    \]
    and $\tilde{R} = \mathcal{E}$.
    According to \Cref{su:equivalent_formulation}, there exists an initial vector $x(1)$ and an input sequence $\{u(k)\}_{k\in\N}$ such that $\{x(k)\}_{k\in\N}$ satisfies all constraints if and only if the system of inequalities \eqref{eq:P-TEG}, in which $R = \tilde{R}\oplus A = A$, is consistent.
    Following \Cref{th:consistency}, this is equivalent to the absence of $\infty$-weight paths in the precedence graph $\graph(M_{\dint{\infty}})$, which coincides with the one represented in \Cref{fi:infinite_precedence_graph}.
    As seen in \Cref{su:precedence_constraints}, although this graph does not contain positive-weight circuits (thus, \eqref{eq:P-TEG} is weakly consistent), there exists an $\infty$-weight path.
    This implies that our problem admits no solution, \ie \eqref{eq:P-TEG} is not consistent.
    Algebraically, this conclusion can be derived by observing that $\Pi_{n^2+1} = \Pi_{5} \neq \Pi_4 = \Pi_{n^2}$, \ie condition (\ref{en:consistent}) in \Cref{th:consistency} is violated; indeed,
    \[
        \Pi_{5} = \begin{bmatrix}
            0 & -\infty\\
            5 & 0
        \end{bmatrix},\
        \Pi_{4} = \begin{bmatrix}
            0 & -\infty\\
            4 & 0
        \end{bmatrix}.
    \]
\end{exmp}


\section{Geometric approach interpretation}\label{se:4}

\subsection{Controlled invariant semimodules}

Given matrices $A\in\Rmax^{n\times n}$ and $B\in\Rmax^{n\times m}$, we say that a semimodule $\X\subseteq \Rmax^n$ is $(A,B)$-invariant (or controlled invariant) if, for all $x(1)\in\X$, there exists a sequence $u(1),u(2),\ldots\in\Rmax^m$ such that $x(2),x(3),\ldots$, obtained through \eqref{eq:max_plus_linear_system}, belong to $\X$.
Suppose that the behavior of a plant to be controlled evolves according to a max-plus linear system \eqref{eq:max_plus_linear_system}, and that we want to impose certain specifications for trajectory $\{x(k)\}_{k\in\N}$, expressed in the form of the inclusion $(\forall k\in\N)\ x(k)\in\mathcal{K}$ where $\mathcal{K}$ is a semimodule.
Then, the maximal $(A,B)$-invariant subsemimodule $\mathcal{K}^*$ of $\mathcal{K}$ is the largest set of vectors $x(1)$ for which there exists a sequence of inputs $u(1),u(2),\ldots$ such that $x(k)\in\mathcal{K}$ for all $k\in\N$ \cite[Lemma 1]{katz2007maxplus}.

It is an open problem to determine, given matrices $A,B$ and a finitely generated semimodule $\mathcal{K}$, the maximal $(A,B)$-invariant subsemimodule $\mathcal{K}^*$ of $\mathcal{K}$.
However, \cite{katz2007maxplus} presents a fixed-point procedure that, if converging in finite time, provides the maximal $(A,B)$-invariant subsemimodule of a given semimodule $\mathcal{K}$.
The procedure is based on the mapping $\phi:2^{\Rmax^n}\rightarrow 2^{\Rmax^n}$ defined by
\[
    \begin{array}{rcl}
        \phi(\X) &=& \X \cap A^{-1}(\X \ominus \Ima B)\\ &=& \X \cap \{x\in\Rmax^{n}\mid (\exists u\in\Rmax^m)\ Ax \oplus Bu\in\X \},
    \end{array} 
\]
where $A^{-1}(\X) = \{x\in\Rmax^n\mid Ax\in\X \}$ and $\X\ominus \mathcal{U} = \{x\in\Rmax^n\mid (\exists u\in\mathcal{U})\ x\oplus u\in\X \}$.
Let us also define $\phi^0(\X) = \X$ and, for all $k\in \N$, $\phi^{k}(\X) = \phi(\phi^{k-1}(\X))$.
If $\mathcal{K}$ is a finitely generated semimodule, then also $\phi(\mathcal{K})$ is, and the elements of the generating matrix $K$, for which $\phi(\mathcal{K}) = \Ima K $, can be computed in finite time.
Moreover, if $\phi^{k+1}(\mathcal{K}) = \phi^k(\mathcal{K})$ for some $k\in\N$, then $\phi^k(\mathcal{K})$ is the maximal $(A,B)$-invariant subsemimodule $\K^*$ of $\mathcal{K}$.
When the sequence does not converge in a finite number of iterations, we know that $\mathcal{K}^*\subseteq \lim_{k\rightarrow +\infty} \phi^k(\mathcal{K})$, but it is unknown whether it is always the case that $\mathcal{K}^* = \lim_{k\rightarrow +\infty} \phi^k(\mathcal{K})$.
For the sake of simplicity, from now on we will indicate $\lim_{k\rightarrow +\infty} \phi^k(\mathcal{K})$ by $\phi^\infty(\mathcal{K})$.


\subsection{Main result}

We will now interpret the results in \Cref{su:solution} in the framework of the geometric approach.
In order to use tools from geometric invariance, we re-write the system dynamics \eqref{eq:scenario_geometric_approach} as in the formulation from \cite{katz2007maxplus}.
Define, for all $k\geq 2$,
\[
    \bar{x}(k) =
    \begin{bmatrix}
    x(k-1) \\ x(k)
    \end{bmatrix},\ 
    \bar{A} = \begin{bmatrix}
        \mathcal{E} & E\\\mathcal{E} & \mathcal{E}
    \end{bmatrix},\
    \bar{B} = \begin{bmatrix}
        \mathcal{E}\\E
    \end{bmatrix},\
    H = \begin{bmatrix}
        C & L\\R & C
    \end{bmatrix}. 
\]
Then, \eqref{eq:scenario_geometric_approach} is equivalently stated as the max-plus linear system
\begin{equation}\label{eq:lifted_max_plus}
    \forall k\geq 2,\quad \bar{x}(k+1) = \bar{A} \bar{x}(k) \oplus \bar{B} u(k),
\end{equation}
subject to the specifications $(\forall k\geq 2)\ \bar{x}(k)\in\mathcal{K}$, where $\mathcal{K}\subseteq \Rmax^{2n}$ is the semimodule
\begin{equation}\label{eq:lifted_semimodule}
    \mathcal{K} = \left\{ \bar{x}\in\Rmax^{2n}\mid \bar{x}\geq H\otimes \bar{x}\right\} = \Ima H^*.
\end{equation}

In the following, we will show that, for the above class of max-plus linear systems and semimodules, strongly polynomial algorithms exist that decide whether $\mathcal{K}^*\cap\R^{2n}$ is non-empty and compute matrix $S$ such that $\Ima S = \mathcal{K}^*$.
The reason for focusing only on real vectors of $\mathcal{K}^*$ is that entries of $x(k)$ equal to $-\infty$ do not have any physical meaning for us, as they would correspond to events that occurred in the infinite past.

Katz showed in \cite[Lemma 6]{katz2007maxplus} that, if $H$ is irreducible, then the geometric control problem can be solved in finite time because the sequence $\{\phi^k(\mathcal{K})\}_{k\in\No}$ converges after a pseudo-polynomial number of iterations.
The following theorems improve this result by extending the analysis to the case in which $H$ is reducible, and by showing that it is possible to completely characterize the convergence of sequence $\{\phi^k(\mathcal{K})\cap\R^{2n}\}_{k\in\No}$ in strongly polynomial time.

\begin{thm}\label{th:geometrical_control}
    Consider the max-plus linear system \eqref{eq:lifted_max_plus} and the semimodule \eqref{eq:lifted_semimodule}.
    Then, $\phi^\infty(\mathcal{K}) = \mathcal{K}^*$.
    Moreover, for all $k\in\No\cup\{\infty\}$, $\phi^k(\K) = \Ima S_{k+2}$, where $S_{k+2}$ is the top-left $2n\times 2n$ block of matrix $M_{\dint{k+2}}^*$.
\end{thm}
\begin{pf}
    By substituting directly the definitions of $\bar{A}$, $\bar{B}$, and $H$ into $\phi(\mathcal{K})$, we get
\begin{align*}
    \phi(\mathcal{K}) &= \mathcal{K} \cap \{\bar{x}\in\Rmax^{2n}\mid (\exists u\in\Rmax^{n})\ \bar{A}\bar{x} \oplus \bar{B}u\in\mathcal{K} \}\\
                      &= \left\{
    \begin{bmatrix}
        x_1\\x_2
    \end{bmatrix}\in\Rmax^{2n}\left| \begin{bmatrix}
        x_1\\x_2
    \end{bmatrix}\geq H\begin{bmatrix}
        x_1\\x_2
\end{bmatrix}\right.\right\} \\
                      &\cap
    \left\{
    \begin{bmatrix}
        x_1\\x_2
    \end{bmatrix}\in\Rmax^{2n}\left|
    (\exists u\in\Rmax^{n})\
    \begin{bmatrix}
        x_2\\u
    \end{bmatrix}  
    \geq H\begin{bmatrix}
        x_2\\u
\end{bmatrix}\right.\right\}\\
                      &=
    \left\{
    \begin{bmatrix}
        x_1\\x_2
    \end{bmatrix}\in\Rmax^{2n}\left|
    \begin{array}{c}
    (\exists x_3\in\Rmax^{n})\\
    \begin{bmatrix}
       x_1\\x_2\\x_3
    \end{bmatrix}  
    \geq 
    \begin{bmatrix}
        C&L&\mathcal{E}\\R&C&L\\\mathcal{E}&R&C
    \end{bmatrix}
    \begin{bmatrix}
        x_1\\x_2\\x_3
\end{bmatrix}
    \end{array} 
\right.\right\},
\end{align*}
where in the last step we renamed $u$ into $x_3$.
By induction, it is immediate to obtain the following expression: 
\[
    \phi^k(\mathcal{K}) = 
    \left\{
    \begin{bmatrix}
        x_1\\x_2
    \end{bmatrix}\in\Rmax^{2n}\left|
    \begin{array}{c}
        (
\exists x_{3},\dots,x_{k+2}\in\Rmax^{n}
    ) \\
    \begin{bsmallmatrix}
       x_1\\x_2\\\vsdots\\x_{k+2}
    \end{bsmallmatrix}
    \geq 
    M_{\dint{k+2}}
    \begin{bsmallmatrix}
        x_1\\x_2\\\vsdots\\x_{k+2}
\end{bsmallmatrix}   
    \end{array} 
    \right.\right\}.
\]
Therefore,
\begin{equation*}
    \phi^\infty(\mathcal{K}) = 
    \left\{
    \begin{bmatrix}
        x_1\\x_2
    \end{bmatrix}\in\Rmax^{2n}\left|
    \begin{array}{c}
     \left(
\exists x_{3},x_4,\ldots\in\Rmax^{n}
    \right)\\
    \begin{bsmallmatrix}
       x_1\\x_2\\x_3\\\vsdots
    \end{bsmallmatrix}
    \geq 
    M_{\dint{\infty}}
    \begin{bsmallmatrix}
        x_1\\x_2\\x_3\\\vsdots
\end{bsmallmatrix}   
    \end{array} 
    \right.\right\},
\end{equation*}
and, since $\phi(\phi^{\infty}(\mathcal{K})) = \phi^\infty(\mathcal{K})$, we have that $\mathcal{K}^* = \phi^\infty(\mathcal{K})$.

It remains to be proven that, for all $k\in\No\cup\{\infty\}$, $\phi^k(\mathcal{K}) = \Ima S_{k+2}$.
The inclusion $\phi^k(\mathcal{K}) \subseteq \Ima S_{k+2}$ is proven by observing that $x_{\dint{k+2}}\geq M_{\dint{k+2}} x_{\dint{k+2}}$ is equivalent to $x_{\dint{k+2}} \geq M_{\dint{k+2}}^* x_{\dint{k+2}}$ and that
\[
    M_{\dint{k+2}}^* x_{\dint{k+2}} \geq \begin{bmatrix}
    S_{k+2} & \mathcal{E}_{2n\times kn}\\\mathcal{E}_{kn\times 2n}&\mathcal{E}_{kn\times kn}
    \end{bmatrix}
     x_{\dint{k+2}},
\]
which implies $\begin{bsmallmatrix}x_1\\x_2\end{bsmallmatrix}\geq S_{k+2}\begin{bsmallmatrix}x_1\\x_2\end{bsmallmatrix}$ or, equivalently, $\begin{bsmallmatrix}x_1\\x_2\end{bsmallmatrix}\in\Ima S_{k+2}^* = \Ima S_{k+2}$.
To prove the inclusion $\phi^k(\mathcal{K})\supseteq \Ima S_{k+2}$, let $\begin{bsmallmatrix}x_1\\x_2\end{bsmallmatrix} \in \Ima S_{k+2}$.
From the definition of the image of a matrix, there exists a vector $\begin{bsmallmatrix}
    u_1\\u_2
\end{bsmallmatrix}\in\Rmax^{2n}$ such that 
\begin{equation}\label{eq:assignment}
    \begin{bmatrix}
        x_1\\x_2
    \end{bmatrix} = S_{k+2} 
    \begin{bmatrix}
        u_1\\u_2
    \end{bmatrix}.
\end{equation}
Now we show that we can find $x_3,x_4,\ldots,x_{k+2}\in\Rmax^{n}$ such that $x_{\dint{k+2}} = [x_1\transposed \ x_2\transposed \ \cdots x_{k+2}\transposed]\transposed$ satisfies $x_{\dint{k+2}}\geq M_{\dint{k+2}}x_{\dint{k+2}}$ or, equivalently, 
\begin{equation}\label{eq:equation}
    x_{\dint{k+2}} =  M_{\dint{k+2}}^*x_{\dint{k+2}}.
\end{equation}
A vector $x_{\dint{k+2}}$ satisfying \eqref{eq:assignment} and \eqref{eq:equation} can be obtained by taking $x_{\dint{k+2}} = M_{\dint{k+2}}^* \tilde{u}$, where $\tilde{u} = \begin{bmatrix}
    u_1\transposed&u_2\transposed&\mathcal{E}_{kn\times 1}\transposed
    \end{bmatrix}\transposed$.
Indeed, \eqref{eq:assignment} comes from the definition of $S_{k+2}$, and \eqref{eq:equation} from 
\[
x_{\dint{k+2}} \! = M_{\dint{k+2}}^*\tilde{u} \! =  M_{\dint{k+2}}^* M_{\dint{k+2}}^*\tilde{u} \! = M_{\dint{k+2}}^* x_{\dint{k+2}}. \qedhere
\]
\end{pf}

For all $k\in\No\cup\{\infty\}$, the following formula for matrix $S_{k+2}$ (\ie the $2n\times 2n$ block in the top-left corner of matrix $M^*_{\dint{k+2}}$) can be obtained using \cite{baccelli1992synchronization}[Lemma 4.101] (see \cite{zorzenon2023switched} and \cite{zorzenon2024infinite} for a proof):
\begin{equation}\label{eq:formula_H}
    S_{k+2} = 
\begin{bmatrix}
    \Pi_{k+1} & \Pi_{k+1} L (\Pi_{k} \oplus \Psi)^*\\
    (\Pi_{k} \oplus \Psi)^* R \Pi_{k+1} & (\Pi_{k} \oplus \Psi)^*
\end{bmatrix},
\end{equation}
where
\[
    \Psi = (R C^* L \oplus C)^*.
\]
Combining \Cref{th:geometrical_control} and the formula \eqref{eq:formula_H}, we get the following result.

\begin{thm}\label{th:geometrical_control_2}
    Consider the max-plus linear system \eqref{eq:lifted_max_plus} and the semimodule \eqref{eq:lifted_semimodule}.
    Then,
\begin{enumerate}
    \item the sequence $\{\phi^k(\mathcal{K})\}_{k\in\No}$ converges in at most $k = n^2-1$ steps and $\mathcal{K}^* \cap\R^{2n} \neq \emptyset$ 
if and only if \eqref{eq:P-TEG} is consistent,\label{en:consistent} 
    \item the sequence $\{\phi^k(\mathcal{K})\cap\R^{2n}\}_{k\in\No}$ does not converge in finitely many steps and $\mathcal{K}^*\cap\R^{2n} = \emptyset$ if and only if \eqref{eq:P-TEG} is weakly consistent but not consistent,\label{en:not_consistent}
    \item the sequence $\{\phi^k(\mathcal{K})\cap\R^{2n}\}_{k\in\No}$ converges in a finite number of steps to the empty set $\emptyset$ if and only if \eqref{eq:P-TEG} is not weakly consistent.\label{en:not_weak}
\end{enumerate}
\end{thm}
\begin{pf}
    Recall, from \Cref{th:geometrical_control}, that $\phi^k(\mathcal{K}) = \Ima S_{k+2}$ for all $k\in\No\cup\{\infty\}$.
\begin{enumerate}
    \item From the equivalence between the geometric control problem and the problem stated in \Cref{se:3}, $\mathcal{K}^*\cap\R^{2n}$ is non-empty if and only if \eqref{eq:P-TEG} is consistent.
        Moreover, if \eqref{eq:P-TEG} is consistent, according to \eqref{eq:formula_H} and \Cref{th:consistency} we have $S_{k+2} = S_{n^2+1}$ for all $k\geq n^2-1$.
        Therefore, $\phi^{k}(\mathcal{K}) = \phi^{n^2-1}(\mathcal{K})$ for all $k\geq n^2-1$.
    \item From \eqref{eq:formula_H} and \Cref{th:consistency}, \eqref{eq:P-TEG} is weakly consistent but not consistent if and only if $\phi^k(\mathcal{K})\cap\R^{2n}\neq \phi^{k+1}(\mathcal{K})\cap\R^{2n}$ for all $k\in\No$.
    \item From \eqref{eq:formula_H} and \Cref{th:consistency}, \eqref{eq:P-TEG} is not weakly consistent if and only if there exists a number $\hat{k}\in\N$ such that $S_{\hat{k}+2}\not\in\Rmax^{2n\times 2n}$.
        To conclude, observe that, because of \Cref{th:gallai}, no real vector belongs to the set $\Ima S_{\hat{k}+2}$ because $S_{\hat{{k}+2}}$ is a star matrix, which implies
    \[
        x \in\Ima S_{\hat{k}+2}\Leftrightarrow x = S_{\hat{k}+2} x\Leftrightarrow x \geq S_{\hat{k}+2} x.\qedhere
    \]
\end{enumerate}
\end{pf}

\begin{rem}
    As shown in \Cref{th:geometrical_control_2}, whenever \eqref{eq:P-TEG} is not consistent, the convergence of $\{\phi^k(\mathcal{K})\}_{k\in\N}$ in finitely many steps is not guaranteed.
However, what is guaranteed is that $\mathcal{K}^*\cap\R^{2n} = \emptyset$.
\end{rem}

\begin{rem}
    In case \ref{en:not_weak} of \Cref{th:geometrical_control_2}, there is a number $\hat{k}\in\N$ such that $\phi^k(\K)\cap\R^{2n} = \emptyset$ for all $k\geq \hat{k}$.
    A pseudo-polynomial upper bound for this number $\hat{k}$ was given in \cite{ZORZENON202219}.
    It is worth mentioning that $\hat{k}$ depends on the magnitude of the entries in $L,C,R$, whereas the algorithms mentioned in \Cref{th:consistency} can be used to determine the convergence of $\{\phi^k(\mathcal{K})\}_{k\in\N}$ in a time that depends only on the dimension $n$ of the matrices.
\end{rem}

\section{Examples}\label{se:5}

\begin{exmp}
Let us take again matrices $L,C,R$ from Example \ref{ex:cont} and define $\bar{A}$, $\bar{B}$, and $\K$ as in \eqref{eq:lifted_max_plus} and \eqref{eq:lifted_semimodule}.
As we have already seen, the system \eqref{eq:P-TEG} is weakly consistent and not consistent.
Therefore, according to \Cref{th:geometrical_control_2}, the maximal $(\bar{A},\bar{B})$-invariant subsemimodule $\K^*$ of $\K$ does not contain any real vector.
Moreover, the sequence $\{\phi^k(\K)\}_{k\in\No}$ does not converge in finitely many steps.
Indeed, $\phi^k(\K) = \Ima S_{k+2}$, and it can be shown (for instance, from the graphical interpretation of the Kleene star) that
\[
    S_{k+2} = 
    \begin{bmatrix}
        0 & -\infty & -\infty &-\infty\\
        1+k & 0 & -1+k &-1\\
        2 & -\infty & 0 &-\infty\\
        2+k & -\infty & k &0
    \end{bmatrix}.
\]
\end{exmp}

\begin{exmp}\label{ex:Katz}
We consider a variation of the transportation network example given in \cite{katz2007maxplus}.
The system evolves according to \eqref{eq:fully_actuated}, where
\[
    A = \begin{bmatrix}
        0 & 17 & -\infty & -\infty\\
        -\infty & 0 & 11 & 9\\
        14 & -\infty & 11 & 9\\
        14 & -\infty & 11 & 0
    \end{bmatrix}.
\]
We consider the single constraint, $\forall k\in\N$, $x_4(k) \geq \ell \otimes x_4(k+1)$, where $\ell\in\R$.
This constraint can be written as \eqref{eq:linear_dual_inequalities} by defining $C = \tilde{R} = \mathcal{E}\in\Rmax^{4\times 4}$, $L_{ij} = \ell$ if $i=j=4$, and $L_{ij} = -\infty$ otherwise.
Defining $R = A\oplus \tilde{R}$ and matrices $\bar{A}$ and $\bar{B}$ as in \eqref{eq:lifted_max_plus}, we want to find, for different values of $\ell$, the maximal $(\bar{A},\bar{B})$-invariant subsemimodule $\K^*$ of $\K = \Ima H^*$, where
\[
    H = \begin{bmatrix}
        C & L\\R & C
    \end{bmatrix} = 
    \begin{bsmallmatrix}
        -\infty & -\infty &-\infty &-\infty &-\infty &-\infty &-\infty &-\infty \\
        -\infty & -\infty &-\infty &-\infty &-\infty &-\infty &-\infty &-\infty \\
        -\infty & -\infty &-\infty &-\infty &-\infty &-\infty &-\infty &-\infty \\
        -\infty & -\infty &-\infty &-\infty &-\infty &-\infty &-\infty &\ell \\
        0 & 17 & -\infty & -\infty&-\infty &-\infty &-\infty &-\infty\\
        -\infty & 0 & 11 & 9&-\infty &-\infty &-\infty &-\infty\\
        14 & -\infty & 11 & 9&-\infty &-\infty &-\infty &-\infty\\
        14 & -\infty & 11 & 0&-\infty &-\infty &-\infty &-\infty
    \end{bsmallmatrix}.
\]
Since $H$ is reducible, we cannot decide the convergence of $\{\phi^k(\K)\}_{k\in\No}$ based on methods developed in \cite{katz2007maxplus}.

\begin{figure}[h]
    \centering
    \begin{tikzpicture}[node distance=2cm and 2cm,>=stealth',bend angle=45,double distance=.5mm,arc/.style={->,>=stealth'},place/.style={circle,fill=gray!15,thick,minimum size=4mm,draw},inner sep=.1mm]

\tiny
\newcommand{\repetitions}{9}
\pgfmathtruncatemacro{\repetitionsmo}{\repetitions-1}
\foreach \z in {1,...,\repetitions}
{
\node [place,draw=white] (n1\z) at (\z*1.01,0) {};
\node [place,draw=white] (n2\z) at (\z*1.01,-1.2) {};
\node [place,draw=white] (n3\z) at (\z*1.01,-2.4) {};
\node [place,draw=white] (n4\z) at (\z*1.01,-3.6) {};
}

\foreach \z in {1,...,\repetitions}
{
\pgfmathtruncatemacro{\za}{\z*4-3}
\pgfmathtruncatemacro{\zb}{\z*4-2}
\pgfmathtruncatemacro{\zc}{\z*4-1}
\pgfmathtruncatemacro{\zd}{\z*4}
\node [place] (n1\z) at (n1\z) {$\za$};
\node [place] (n2\z) at (n2\z) {$\zb$};
\node [place] (n3\z) at (n3\z) {$\zc$};
\node [place] (n4\z) at (n4\z) {$\zd$};
}

\foreach[evaluate=\z as \zz using int(\z+1)] \z in {1,...,\repetitionsmo}
{
\draw [arc] (n1\z) to node[auto] {} (n1\zz);
\draw [arc] (n2\z) to node[auto] {} (n2\zz);
\draw [arc] (n3\z) to node[auto,pos=.2] {$11$} (n3\zz);
\draw [arc] (n4\z) to node[auto] {} (n4\zz);

\draw [arc] (n4\z) to node[auto,swap,pos=.25] {$9$} (n3\zz);
\draw [arc] (n3\z) to node[auto,swap,pos=.15] {$11$} (n4\zz);
\draw [arc] (n1\z) to node[auto,pos=.1] {$14$} (n3\zz);
\draw [arc] (n3\z) to node[auto,pos=.25] {$11$} (n2\zz);

\ifnum\z=1
\else
\draw [arc] (n4\z) to node[auto,pos=.12] {$9$} (n2\zz);
\fi
\ifnum\z=2
\else
\draw [arc] (n2\z) to node[auto,swap,pos=.6] {$17$} (n1\zz);
\fi
\ifnum\z=3
\else
\draw [arc] (n1\z) to node[auto,swap,pos=.4] {$14$} (n4\zz);
\fi
}
\foreach[evaluate=\z as \zz using int(\z-1)] \z in {5,...,\repetitions}
{
\draw [arc] (n4\z) to [bend left] node[auto] {$\ell$} (n4\zz);
}

\draw [arc,dashed,thick] (n41) to node[auto,pos=.12] {$9$} (n22);
\draw [arc,dashed,thick] (n22) to node[auto,swap,pos=.6] {$17$} (n13);
\draw [arc,dashed,thick] (n13) to node[auto,swap,pos=.4] {$14$} (n44);
\draw [arc,dashed,thick] (n44) to [bend left] node[auto] {$\ell$} (n43);
\draw [arc,dashed,thick] (n43) to [bend left] node[auto] {$\ell$} (n42);
\draw [arc,dashed,thick] (n42) to [bend left] node[auto] {$\ell$} (n41);


\draw [glow=myred,draw=none] (n41) to (n22) to (n13) to (n34) to (n25) to (n16) to (n47) to [bend left] (n46);
\draw [glow=myred,draw=none] (n46) to [bend left] (n45);
\draw [glow=myred,draw=none] (n45) to [bend left] (n44);
\draw [glow=myred,draw=none] (n44) to [bend left] (n43);
\draw [glow=myred,draw=none] (n43) to [bend left] (n42);
\draw [glow=myred,draw=none] (n42) to [bend left] (n41);

\fill [white] (\repetitions-.5,.5) rectangle (\repetitions+.5,-4.1);
\node at (\repetitions,-1.7) {\textbf{\dots}};

\draw [arc] (1,-4.1) to node[pos=.95,above] {$\N$} (\repetitions+.5,-4.1);

\end{tikzpicture}%
    \caption{Precedence graph $\graph(M_{\dint{\infty}})$ for Example \ref{ex:Katz}. Where not indicated, the weight of arcs is $0$.}\label{fi:infinite_precedence_graph_Katz}
\end{figure}

For $\ell = -14$, we get 
\[
\Pi_{16} = \Pi_{17} = \begin{bmatrix}
    0 & -\infty  &-\infty  &-\infty\\
    -\infty     &0&  -\infty & -\infty\\
    -\infty & -\infty  &   0 & -\infty\\
    0   &  3  &   0  &   0
\end{bmatrix}\in\Rmax^{4\times 4},
\]
which implies that \eqref{eq:P-TEG} is consistent.
Hence, according to \Cref{th:geometrical_control_2}, $\{\phi^k(\K)\}_{k\in\No}$ converges in at most $15$ steps (in this example, convergence is reached in only $2$ steps) and $\K^*\cap\R^8\neq\emptyset$.
On the other hand, taking $\ell>-14$ makes the system \eqref{eq:P-TEG} not weakly consistent.
For instance, for $\ell = -13$, the graph $\graph(M_{\dint{\infty}})$ contains the positive-weight circuit formed by the dashed arcs in \Cref{fi:infinite_precedence_graph_Katz}; the convergence of $\{\phi^k(\K)\cap\R^8\}_{k\in\No}$ to $\emptyset$ occurs for $k = 2$.
Interestingly, for values of $\ell>-14$ closer to $-14$, the convergence of $\{\phi^k(\K)\cap\R^8\}_{k\in\No}$ becomes slower, because the length of the shortest positive-weight circuit increases.
Numerically, we verified that, for $\ell=-13.5$, convergence occurs in $5$ steps (see the arcs highlighted in red in \Cref{fi:infinite_precedence_graph_Katz}); for $\ell=-13.9$ in $20$ steps; for $\ell = -13.999$ in $2000$ steps.
This shows the advantage of the strongly polynomial-time algorithm defined in \cite{ZORZENON202219} for checking weak consistency, since its computational time does not depend on the magnitude of elements in matrix $H$.
\end{exmp}

\section{Conclusions}\label{se:conclusions}

The results presented in this paper expand the class of geometric control problems in max-plus linear systems for which a solution can be computed in finite time.
Namely, we have shown that, for fully actuated max-plus linear systems subject to precedence semimodules, it is possible to compute in strongly polynomial time the matrix whose image is the maximal controlled invariant subsemimodule.

For the class of systems and semimodules investigated in this paper, another interesting property can be verified: all controlled invariants are of feedback type, \ie whenever the maximal controlled invariant subsemimodule is nonempty, there is a static max-plus linear feedback that keeps the state inside this set.
The proof of this fact will be given in future work.

\bibliography{references}

\begin{thebibliography}{17}
\providecommand{\natexlab}[1]{#1}
\providecommand{\url}[1]{\texttt{#1}}
\providecommand{\urlprefix}{URL }
\expandafter\ifx\csname urlstyle\endcsname\relax
  \providecommand{\doi}[1]{doi:\discretionary{}{}{}#1}\else
  \providecommand{\doi}{doi:\discretionary{}{}{}\begingroup
  \urlstyle{rm}\Url}\fi

\bibitem[{Animobono et~al.(2023)Animobono, Scaradozzi, Zattoni, Perdon, and
  Conte}]{9832498}
Animobono, D., Scaradozzi, D., Zattoni, E., Perdon, A.M., and Conte, G. (2023).
\newblock The model matching problem for max-plus linear systems: A geometric
  approach.
\newblock \emph{IEEE Transactions on Automatic Control}, 68(6), 3581--3587.
\newblock \doi{10.1109/TAC.2022.3191362}.

\bibitem[{Baccelli et~al.(1992)Baccelli, Cohen, Olsder, and
  Quadrat}]{baccelli1992synchronization}
Baccelli, F., Cohen, G., Olsder, G.J., and Quadrat, J.P. (1992).
\newblock \emph{{Synchronization and linearity: an algebra for discrete event
  systems}}.
\newblock John Wiley \& Sons Ltd.

\bibitem[{Basile and Marro(1991)}]{basile2002controlled}
Basile, G. and Marro, G. (1991).
\newblock \emph{Controlled and conditioned invariants in linear system theory}.
\newblock Prentice Hall.

\bibitem[{Conte and Perdon(1995)}]{conte1995disturbance}
Conte, G. and Perdon, A.M. (1995).
\newblock The disturbance decoupling problem for systems over a ring.
\newblock \emph{SIAM Journal on Control and Optimization}, 33(3), 750--764.

\bibitem[{Cormen et~al.(2022)Cormen, Leiserson, Rivest, and
  Stein}]{cormen2022introduction}
Cormen, T.H., Leiserson, C.E., Rivest, R.L., and Stein, C. (2022).
\newblock \emph{Introduction to algorithms}.
\newblock MIT press.

\bibitem[{Di~Loreto et~al.(2010)Di~Loreto, Gaubert, Katz, and
  Loiseau}]{di2010duality}
Di~Loreto, M., Gaubert, S., Katz, R.D., and Loiseau, J.J. (2010).
\newblock Duality between invariant spaces for max-plus linear discrete event
  systems.
\newblock \emph{SIAM Journal on Control and Optimization}, 48(8), 5606--5628.

\bibitem[{Gallai(1958)}]{gallai1958maximum}
Gallai, T. (1958).
\newblock Maximum-minimum {S}{\"a}tze {\"u}ber {G}raphen.
\newblock \emph{Acta Mathematica Hungarica}, 9(3-4), 395--434.

\bibitem[{Gonçalves et~al.(2016)Gonçalves, Maia, and
  Hardouin}]{goncalves2016onthe}
Gonçalves, V.M., Maia, C.A., and Hardouin, L. (2016).
\newblock On the steady-state control of timed event graphs with firing date
  constraints.
\newblock \emph{IEEE Transactions on Automatic Control}, 61(8), 2187--2202.
\newblock \doi{10.1109/TAC.2015.2481798}.

\bibitem[{Hardouin et~al.(2018)Hardouin, Cottenceau, Shang, and
  Raisch}]{hardouin2018control}
Hardouin, L., Cottenceau, B., Shang, Y., and Raisch, J. (2018).
\newblock Control and state estimation for max-plus linear systems.
\newblock \emph{Foundations and Trends® in Systems and Control}, 6(1), 1--116.
\newblock \doi{10.1561/2600000013}.

\bibitem[{Katz(2007)}]{katz2007maxplus}
Katz, R.D. (2007).
\newblock Max-plus {$(A,B)$}-invariant spaces and control of timed
  discrete-event systems.
\newblock \emph{IEEE Transactions on Automatic Control}, 52(2), 229--241.
\newblock \doi{10.1109/TAC.2006.890478}.

\bibitem[{Khansa et~al.(1996)Khansa, Denat, and
  Collart-Dutilleul}]{khansa1996a}
Khansa, W., Denat, J.P., and Collart-Dutilleul, S. (1996).
\newblock Structural analysis of {P}-time {P}etri nets.
\newblock In \emph{CESA '96 IMACS Multiconférence, Lille France July 9-12},
  127--136.

\bibitem[{Maia et~al.(2011)Maia, Andrade, and Hardouin}]{maia2011control}
Maia, C.A., Andrade, C., and Hardouin, L. (2011).
\newblock On the control of max-plus linear system subject to state
  restriction.
\newblock \emph{Automatica}, 47(5), 988--992.

\bibitem[{Martinez et~al.(2022)Martinez, Kara, Abdesselam, and
  Loiseau}]{martinez2022systems}
Martinez, C., Kara, R., Abdesselam, A.N., and Loiseau, J.J. (2022).
\newblock Systems synchronisation in max-plus algebra: a controlled invariance
  perspective in memoriam {\'e}douard wagneur.
\newblock \emph{IFAC-PapersOnLine}, 55(40), 1--6.

\bibitem[{Wonham(1974)}]{wonham1974linear}
Wonham, W.M. (1974).
\newblock \emph{Linear multivariable control}, volume 101.
\newblock Springer.

\bibitem[{Zorzenon et~al.(2022)Zorzenon, Balun, and Raisch}]{ZORZENON202219}
Zorzenon, D., Balun, J., and Raisch, J. (2022).
\newblock Weak consistency of {P}-time event graphs.
\newblock \emph{IFAC-PapersOnLine}, 55(40), 19--24.
\newblock \doi{10.1016/j.ifacol.2023.01.042}.
\newblock 1st IFAC Workshop on Control of Complex Systems COSY 2022.

\bibitem[{Zorzenon et~al.(2024)Zorzenon, Komenda, and
  Raisch}]{zorzenon2023switched}
Zorzenon, D., Komenda, J., and Raisch, J. (2024).
\newblock Switched max-plus linear-dual inequalities: cycle time analysis and
  applications.
\newblock \emph{Discrete Event Dynamic Systems}, 34(1), 199--250.

\bibitem[{Zorzenon and Raisch(2025)}]{zorzenon2024infinite}
Zorzenon, D. and Raisch, J. (2025).
\newblock Infinite precedence graphs for consistency verification in {P}-time
  event graphs.
\newblock Submitted for publication.

\end{thebibliography}

\end{document}